\documentclass[12pt]{amsart}
\usepackage{amsfonts,latexsym,amsthm,amssymb,graphicx}
\usepackage[all]{xy}

\newtheorem{theorem}{Theorem}[section]

{\theoremstyle{definition} \newtheorem{defin}[theorem]{Definition}}
{\theoremstyle{remark} 
\newtheorem{example}[theorem]{Example}}

\newcommand{\Abb}{{\mathbb{A}}}

\newcommand{\Lbb}{{\mathbb{L}}}
\newcommand{\Pbb}{{\mathbb{P}}}
\newcommand{\Qbb}{{\mathbb{Q}}}
\newcommand{\Zbb}{{\mathbb{Z}}}

\DeclareMathOperator{\dCan}{dCan}

\DeclareMathOperator{\Ob}{Ob}
\DeclareMathOperator{\ord}{ord}

\newcommand{\one}{1\hskip-3.5pt1}
\newcommand{\csm}{c_{\text{\rm SM}}}
\newcommand{\integ}[3][X]{\int_{#3}\one(#2)\,d\mathfrak c_#1}

\newcommand{\cA}{{\mathcal A}}
\newcommand{\cC}{{\mathcal C}}
\newcommand{\cD}{{\mathcal D}}
\newcommand{\cK}{{\mathcal K}}
\newcommand{\cS}{{\mathcal S}}
\newcommand{\cV}{{\mathcal V}}
\newcommand{\cX}{{\mathcal X}}
\newcommand{\cY}{{\mathcal Y}}

\newcommand{\id}{{\text{id}}}

\title[Celestial integration]{Celestial integration, stringy
invariants, and Chern-Schwartz-MacPherson classes}
\author{Paolo Aluffi}
\address{Max-Planck-Institut f\"ur Mathematik, Bonn, Deutschland}
\address{Math Dept, Florida State University, Tallahassee, Florida, 
U.S.A.}
\email{aluffi@math.fsu.edu}

\begin{document}

\begin{abstract}
We introduce a formal integral on the system of varieties mapping
properly and birationally to a given one, with value in an associated
Chow group. Applications include comparisons of Chern numbers of
birational varieties, new birational invariants, `stringy' Chern
classes, and a `celestial' zeta function specializing to the
topological zeta function.

In its simplest manifestation, the integral gives a new expression for
Chern-Schwartz-MacPherson classes of possibly singular varieties,
placing them into a context in which a `change of variable' formula
holds.

The formalism has points of contact with motivic integration.
\end{abstract}

\maketitle


\section{Introduction}\label{intro}

\subsection{}
In this note I review the notion of {\em celestial integration,\/}
and sketch a few applications; for proofs and further details, the
reader is addressed to \cite{math.AG/0407150}. I am very grateful to
Jean-Paul and Cidinha for organizing the {\em VIII\`eme Rencontre
Internationale de {\em S\~ao Carlos\/} sur les singularit\'es
re\'elles et complexes au CIRM,\/} relocating for the occasion the
idyllic surroundings of {\em S\~ao Carlos, Brasil\/} to the idyllic
surroundings of {\em Luminy, France.\/} Perfect weather, exceptionally
interesting talks, and spirited conversations made the conference a
complete success. What follows was the subject of my lecture at the
S\~ao Carlos/Luminy meeting, and preserves (for better or worse) the
informal nature of a seminar talk.

I thank the Max-Planck-Institut f\"ur Mathematik in Bonn, where much
of this material was conceived and where this note was written.

\subsection{}
Summary: to a variety $X$ I will associate a large group $A_*\cC_X$
(containing the Chow group $A_*X$ of $X$); for certain data $\cD$,
$\cS$ (arising, for example, from a divisor $D$ and a constructible
subset $S$ of $X$) I will define a distinguished element
$$\integ{\cD}{\cS} \in A_*\cC_X\quad.$$ 
These are the {\em celestial integrals\/} in the title; they are not
defined as integrals, but satisfy formal properties justifying the
terminology.  The qualifier {\em celestial\/} is meant to evoke the
fact that the {\em modification systems\/} on which this operation is
defined are close relatives of Hironaka's {\em vo\^ute \'etoil\'ee.\/}

Applications of this construction:
\begin{itemize}
\item Comparison of Chern classes of birational varieties;
\item Birational invariants;
\item `Celestial' zeta functions;
\item Invariants of singular varieties ({\em `Stringy Chern classes'\/});
\item Relations with the theory of Chern-Schwartz-MacPherson classes.
\end{itemize}

\subsection{}
Digression: {\bf motivic integration}. This is technically not
necessary for the rest of the talk, but useful nonetheless as
`inspiration' for the main construction.

The {\em Grothendieck group of varieties\/} is the free abelian group
on symbols $[X]$, where $X$ is a complex algebraic variety up to
isomorphism, modulo relations
$$[X]=[Y]+[X\setminus Y]$$
for each closed subvariety $Y\subset X$. This group may be made into a
ring by setting $[X]\cdot [Y]=[X\times Y]$.

\begin{example}
The class of a point is the identity for this multiplication. The
class $[\Abb^1]$ is denoted $\Lbb$; thus
$$[\Pbb^n]=[\Abb^0]+[\Abb^1]+\dots+[\Abb^n]=\frac{\Lbb^{n+1}-1}{\Lbb-1}
\quad.$$
\end{example}
As Eduard Looijenga writes (\cite{MR2003k:14010}, p.~269), this
ring---or rather its localization at $\Lbb$---is {\em interesting,
big, and hard to grasp.\/} In practice, it is necessary to further
tweak this notion, by a suitable completion with respect to a
dimension filtration; I will glibly ignore such important `details'.

Mapping $X$ to its class in this ring gives a {\em universal Euler
characteristic:\/} anything satisfying the basic relations---e.g.,
topological Euler characteristic, Hodge polynomials and
structures,\dots, must factor through this map. This is motivation to
`compute' $[X]$ for given $X$.

Through {\em motivic integration,\/} one can determine an element 
$$\int_\cS\Lbb^{-\ord D} d\mu$$
in the completed Grothendieck ring of varieties, from the information
of a divisor $D$ and a constructible subset $\cS$ of the arc space
$\cX$ of $X$.  If $X$ is nonsingular, then choosing $D=0$, $\cS=\cX$
gives $\int_\cX \Lbb^{0} d\mu = [X]$.

Motivic integration was defined and developed by Maxim Kontsevich, Jan
Denef, and Fran\c cois Loeser, and it is of course much deeper than
this brief summary can begin to suggest. There are many good surveys
of this material, for example: \cite{math.AG/9911179},
\cite{MR1905328}, \cite{MR2003k:14010}, \cite{math.AG/0401374}. For an
explanation of what makes motivic integration {\em motivic,\/} see the
appendix in \cite{math.AG/9911179}.

The definition of the integral $\int_\cS\Lbb^{-\ord D} d\mu$ relies on
the study of the {\em arc spaces\/} of a variety; $d\mu$ is a measure
on this space, with value in the completed Grothendieck ring; the
integral is an honest integral with respect to this measure, and as
such it satisfies a {\em change of variable\/} formula: for example,
if $\pi: Y\to X$ is proper and birational, then
$$\int_\cX \Lbb^{-\ord D} d\mu = \int_\cY \Lbb^{-\ord(\pi^{-1} D+K_\pi)} 
d\mu\quad,$$
where $K_\pi$ denotes the relative canonical sheaf.

This formula is at the root of spectacular applications of motivic
integration.  For example, suppose $X$, $Y$ are {\em birational\/}
nonsingular, complete Calabi-Yau varieties; resolve a birational
morphism between them:
$$\xymatrix{
& V \ar[dl]_{\pi_1} \ar[dr]^{\pi_2} \\
X \ar@{-->}[rr] & & Y
}$$
with $V$ nonsingular and $\pi_1$, $\pi_2$ proper and birational. Then
the Calabi-Yau condition implies $K_{\pi_1}=K_{\pi_2}$; denoting this
by $K$ gives
$$[X]=\int_\cX \Lbb^{0} \,d\mu = \int_\cV\Lbb^{-\ord K} \,d\mu=
\int_\cY \Lbb^{0} \, d\mu=[Y]\quad.$$
Hence: such varieties must have the same topological Euler
characteristic, Betti numbers, Hodge polynomials, etc.

\subsection{}
\underline{Motiv}ic integr\underline{ation} only serves as {\em
motivation\/} for the rest of this lecture, or maybe more correctly as
a motivating analogy.  The basic relation in the Grothendieck ring
holds (in a suitable sense) for {\em Chern-Schwartz-MacPherson\/}
classes; `hence' there should be a `motivic' theory of such classes:
it should be possible to deal with the classes within the framework of
an integration theory, satisfying a suitable change-of-variable
formula; one should be able to play tricks such as the application
sketched above at the level of Chern classes.

This is the guiding theme in what follows.


\section{Modification systems}

\subsection{}
The task is to define an `integral' carrying information about Chern
classes.  Taking at heart the lesson learned in motivic integration,
we should start by defining an appropriate context in which this
integral may take its value.

Let $X$ be a variety over an algebraically closed field of
characteristic zero (the precise requirement is that embedded
resolution \`a la Hironaka should work).

\begin{defin}
I will denote by $\cC_X$ the category of proper birational maps
$$\xymatrix{
V_\pi \ar[d]^\pi \\
X
}$$
with morphisms given in the obvious way by commutative triangles
$$\xymatrix{
V_{\pi_1} \ar[rr]^\alpha\ar[dr]_{\pi_1} & & V_{\pi_2} \ar[ld]^{\pi_2}\\
& X
}$$
with $\alpha$ proper and birational.
\end{defin}

Proper birational maps are often called {\em modifications,\/} and the
natural way to think of $\cC_X$ is as an inverse system, so it seems
appropriate to call this category the {\em modification system\/} of
$X$.  Also, it is useful to take this notion up to the following
equivalence relation: say that $\cC_X$ and $\cC_Y$ are equivalent if
there exists objects in $\cC_X$ and $\cC_Y$ with a common source. For
example, if $X$ and $Y$ are birational and complete then their
modification systems are equivalent in this sense.

\subsection{}
I will (usually) denote by $V_\pi$ the source of the object $\pi$ of
$\cC_X$.  It is hard to resist the temptation to think of the object
$\pi$ really in terms of its corresponding $V_\pi$, and of $\cC_X$ as
a system of varieties birational to $X$.

Morally I would like to take the inverse limit of this system, and
define ordinary data such as divisors, Chow group, etc.~for the
resulting {\em pro\/}variety.  In practice, it is more straightforward
to simply define these data as appropriate limits of the corresponding
data on the individual $V_\pi$'s.  For example, denote by $A_*V_\pi$
the Chow group of $V_\pi$, with rational coefficients; then
$$\cA_X:=\{A_*V_\pi\, |\, \pi\in \Ob(\cC_X)\}$$
is an inverse system of abelian groups under proper push-forward.

\begin{defin}
The {\em Chow group\/} of $\cC_X$ is the inverse limit of this system:
$$A_*\cC_X:=\varprojlim \cA_X\quad.$$
\end{defin}

Thus, an element $a\in A_*\cC_X$ consists of the data of a class
$(a)_{\id}$ in the Chow group of $X$ and of compatible lifts
$(a)_{\pi}$ for all $\pi\in\Ob(\cC_X)$. I call $(a)_{\pi}$ the {\em
$\pi$-manifestation\/} of $a$.

Note that (if e.g., $X$ is nonsingular) any class $\alpha\in A_*X$ 
determines a `silly' class $a\in A_*\cC_X$: just set $(a)_\pi:=
\pi^* \alpha$. One intriguing (to me, at
least) consequence of the construction given in this paper is that
certain classes on $X$ have {\em other,\/} more interesting, lifts to
$A_*\cC_X$.  These lifts call for rational coefficients, hence the
need for rational coefficients in the definition of $A_*\cC_X$.

Equivalent modification systems have isomorphic Chow groups. 

\subsection{}\label{divisors}
Other standard notions may be defined similarly. Divisors and
constructible sets of sources $V_\pi$ are organized by {\em direct\/}
systems, under pull-backs; the corresponding notions for a
modification system are defined as direct limits of these systems.

For example, a divisor $\cD$ of $\cC_X$ is represented by a pair
$(\pi, D_\pi)$ with $D_\pi$ a divisor of $V_\pi$, and where pairs
$(\pi, D_\pi)$, $(\pi\circ\alpha, D_{\pi\circ\alpha})$ are identified
whenever $\alpha: V_{\pi\circ \alpha} \to V_\pi$ is a proper
birational map and $D_{\pi\circ\alpha}=\alpha^{-1}(D_\pi)$:
$$\xymatrix{
V_{\pi\circ\alpha} \ar[r]^\alpha & V_\pi \ar[r]^\pi & X\quad.
}$$

An obvious way to get a divisor $\cD$ is by pulling back a divisor of
$X$ through the whole system; but note that there are {\em many
other\/} divisors: for example, every subscheme $S$ of $X$ determines
a divisor of $\cC_X$ (represented by the exceptional divisor in the
blow-up of $X$ along $S$).  As a bonus, equivalent modification
systems have the same divisors, while birational varieties don't.

\subsection{}
The story is entirely analogous for {\em constructible subsets\/} of a
modification system. The `obvious' such object is determined by a
constructible (for example, closed) subset of $X$, by taking inverse
images through the system. While this is our main example, one should
keep in mind that the notion is considerably more general.

Of course, equivalent systems have the same constructible subsets. For 
example, if $V$ maps properly and birationally to both $X$ and $Y$:
$$\xymatrix{
& V \ar[dl]_{\pi_X} \ar[dr]^{\pi_Y} \\
X \ar@{-->}[rr] & & Y
}$$
then $(\pi_X, V)$ determines the same subset of $\cC_X$ as $(\id,X)$
and the same subset of $\cC_Y$ as $(\id,Y)$; abusing language I may
denote this object by $\cC_X$ or $\cC_Y$ according to the context, but
the reader should keep in mind that these constructible subsets of
different systems may be identified.

Details about all these notions, and natural definitions (such as sums
of divisors, or unions of constructible subsets) are left to the 
interested reader, and may be found in \cite{math.AG/0407150}.


\section{Celestial integrals}\label{main}

\subsection{}\label{deff}
The main result of this note is that for a variety $X$, a divisor
$\cD$ of $\cC_X$, and a constructible subset $\cS$, there is an
element (the `celestial integral' of $\cD$ over $\cS$, in $\cC_X$)
$$\integ{\cD}{\cS}$$
of the Chow group $A_*\cC_X$, satisfying interesting properties.

The actual definition of this element is uninspiring; I'll give it at
the end of the paper for the sake of completeness. The properties
satisfied by this `integral' are more important, so they get the honor
of prime time.

Of course the celestial integral is additive with respect to disjoint unions of
constructible subsets, as should be expected from an integral. What
makes the notion interesting is that it computes interesting objects
for suitable choices of the input data, and that it satisfies a
change-of-variable formula (again, as should be expected from an
integral!).

More explicitly:
\begin{theorem}\label{props}
\begin{enumerate}
\item\label{norm} {\em (Normalization)\/} Assume $X$ is nonsingular. 
If $\cS$ is represented by $(\id, S)$ with $S\subset X$ a nonsingular
subvariety, then
$$\left(\integ{0}{\cS}\right)_{\id}=c(TS)\cap [S]\quad,$$
the {\em total homology Chern class\/} of $S$, viewed as an element of
$A_*X$.
\item\label{cov} {\em (Change-of-variables)\/} If $\rho: Y \to X$ is
proper and birational, then
$$\integ{\cD}{\cS}=\integ[Y]{\cD+K_\rho}{\cS}\quad,$$
where $K_\rho$ denotes the relative canonical divisor of $\rho$.
\end{enumerate}
\end{theorem}
The normalization property (\ref{norm}) is self-explanatory; it will
serve as a point of depart for extensions to possibly {\em singular\/}
subsets $S$ of a nonsingular variety, in \S\S\ref{stringin}
and~\ref{CSM}. By contrast, (\ref{cov}) needs an immediate
clarification.

\subsection{}
First of all, note that under the given hypotheses we have that
$\cC_X$ and $\cC_Y$ are equivalent systems; thus we may indeed treat
$\cS$ and $\cD$ as data belonging to either.

Secondly, I have to clarify what I mean by {\em the relative
canonical divisor\/} $K_\rho$ of $\rho: Y\to X$. It in fact turns out
that there are more than one sensible such notions, according to the
context. In the simplest case, when $X$ and $Y$ are nonsingular,
$K_\rho$ is the divisor corresponding to the determinant of the
differential $d\rho: TY \to\rho^* TX$. In the general case there is a
choice as to what is the correct generalization of $TX$; I'll come
back to this point in \S\ref{diffi}.

Further, in general the scheme corresponding to the vanishing of the
determinant may not be locally principal. This is not a problem in our
context, however, since every subscheme of every variety in the
modification system $\cC_X$ determines a divisor of the system, as
observed in \S\ref{divisors}.


\section{Sketch of applications}\label{applic}

\subsection{Invariance of Chern classes}
Exactly as in the case of motivic integration, the change-of-variable
formula yields an invariance statement for celestial integration
across birational morphisms preserving the canonical class.

Two varieties $X$, $Y$ (nonsingular, for simplicity) are {\em
$K$-equivalent\/} if their modification systems are equivalent, and
the canonical divisors $K_X$, $K_Y$ agree in the system(s): that is,
if the pull-backs of $K_X$, $K_Y$ to a common source agree:
$$\xymatrix{
& V \ar[dl]_{\pi_X} \ar[dr]^{\pi_Y} \\
X \ar@{-->}[rr] & & Y
}$$
$$\qquad\pi_X^* K_X = \pi_Y^* K_Y\quad.$$

In this situation, letting $K=K_{\pi_X}=K_{\pi_Y}$ and applying the 
change-of-variable formula (\ref{cov}) gives:
$$\integ{\cD}{\cS}=\integ[V]{\cD+K}{\cS} =\integ[Y]{\cD}{\cS}\quad.$$
Therefore:
\begin{theorem}
Celestial integrals on $K$-equivalent varieties agree as elements of
the (common) Chow group of the corresponding modification systems.
\end{theorem}

\subsection{}
For example, applying this observation with $\cD=0$, $\cS=\cC_X$
($=\cC_Y$ on $Y$) and using (\ref{norm}) from \S\ref{main} shows that
$c(TX)\cap [X]$ and $c(TY)\cap [Y]$ are manifestations (in $A_*X$,
$A_*Y$, respectively) of {\em the same\/} class in the Chow group of
the modification system.

This recovers the fact (known, for example, through motivic
integration) that the Euler characteristics of $K$-equivalent
varieties must agree.  More generally, it shows (via a simple
application of the projection formula) that all numbers
$$c_1^i\cdot c_{n-i}$$
with $n=\dim X=\dim Y$ must agree for $K$-equivalent varieties.

Incidentally, these numbers must therefore be invariant through
classical flops, and hence (as shown by Burt Totaro,
\cite{MR2001g:58037}) they must factor through the complex elliptic
genus. It is a pleasant exercise to verify this fact directly: these
numbers can be assembled into a genus (which I would like to call the
{\em cuspidal\/} genus, for reasons which will likely be apparent to
many readers), corresponding to the {\em characteristic\/}
$e^{xT}(1+xU)$; it is straightforward to check directly that the
cuspidal genus factors through the complex elliptic genus.

The invariance of Chern numbers mentioned above is of course only a
very particular case of similar results accessible through celestial
integration. Every choice of a divisor and a constructible subset
yields an analogous (but, unfortunately, usually much less
transparent) invariance statement.

\subsection{Birational invariants}
As another application, we can extract new birational invariants from
the integral. For example, let
$$\dCan(X):=\{\deg \integ{\cK}{\cC_X}\, |\, \text{$\cK$ effective
  canonical divisor of $X$}\} \subset \Zbb\quad.$$

\begin{theorem}
If $X$, $Y$ are birational complete varieties, then $\dCan(X)=\dCan(Y)$.
\end{theorem}

For example: if $X$ is birational to a Calabi-Yau variety $Y$, then
$\dCan(X)=\{\chi(Y)\}$.

An interesting question is whether an analogous (and nontrivial)
invariant can be defined for varieties without effective canonical
divisors. A few natural candidates for such invariants, involving
negative representatives, must be ruled out at the moment because of a
sticky technical obstacle to the definition of celestial integrals for
noneffective divisors (see \S\ref{defin}).

\subsection{Zeta functions}
In analogy with motivic integration, zeta functions can be concocted
from motivic integrals. For example, for a divisor $D$ of $X$ (say
defined by $f=0$) set
$$Z(D,m):=\integ{m\cD}{\cC_X}\quad,$$
a series in the variable $m$, with coefficients in $\cA_*\cC_X$;
here $\cD$ is the divisor in $\cC_X$ determined by $D$. Then

\begin{theorem}
The degree of $Z(D,m)$ equals the topological zeta function of $f$.
\end{theorem}

This connection makes it possible to formulate analogs of the {\em
monodromy conjecture\/} (see for example \cite{math.AG/0401374},
\S6.8) for celestial zeta functions. I hope that the celestial
viewpoint will add something to the circle of ideas surrounding zeta
functions. For example, conceivably the relationship between celestial
integration and the theory of Chern-Schwartz-MacPherson classes
(\S\ref{CSM}) may give a tool to compute local contributions to the
zeta function of a hypersurface in terms of the Segre class of its
singularity subscheme.


\section{Stringy invariants}\label{stringin}

\subsection{}\label{stringy}
If $X$ has sufficiently mild singularities, there is a notion of {\em
stringy Euler characteristic\/} of $X$, introduced by Batyrev. For
example, in the particular case in which $X$ admits a {\em crepant\/}
resolution $V$, the stringy Euler characteristic of $X$ may be defined
to be the ordinary Euler characteristic $\chi(V)$ of $V$: remarkably,
this turns out to be independent of the chosen crepant resolution.

Celestial integration extends this notion to a whole class in $A_*X$. 
By the normalization property (Theorem~\ref{props} (\ref{norm})),
$$\left(\integ{0}{\cC_X}\right)_{\id}=c(TX)\cap [X]$$ 
if $X$ is nonsingular; but the expression on the left-hand-side
defines an element of $A_*X$ {\em even if $X$ is singular\/} (in fact,
the celestial integral defines this expression together with
distinguished lifts to all varieties mapping to $X$). If $X$ admits a
crepant resolution $\pi: V \to X$, it is easy to check that this
definition produces the push-forward $\pi_* c(TV)\cap [V]$. By the
Poincar\'e-Hopf theorem, therefore, the degree of this class recovers
the stringy Euler characteristic of $X$ in this case.

This in fact holds for any $X$ for which the stringy Euler
characteristic is defined, justifying the following:
\begin{defin}
The {\em stringy Chern class\/} of $X$ is the identity manifestation
$$\left(\integ{0}{\cC_X}\right)_{\id}\quad.$$
\end{defin}

Coincidentally, a notion of stringy Chern class was produced
simultaneously as the one presented above, by Tommaso de Fernex,
Ernesto Lupercio, Thomas Nevins, and Bernardo Uribe (in fact, the
preprint \cite{math.AG/0407314} appeared on the arXiv during the
S\~ao~Carlos/Luminy conference!). While the approaches to the two
notions differ somewhat, the two stringy classes agree.

\subsection{}\label{diffi}
There is a subtlety here, which I can only touch upon in this
note. For singular $X$, the notion of celestial integral depends on
the choice of a good notion of {\em relative canonical divisor\/.} The
`usual' notion is constructed starting from the {\em double dual\/}
$\omega_X$ of the K\"ahler differentials $\Omega^{\dim X}_X$ of $X$;
this `$\omega$ flavor' of the celestial integral is what leads to the
stringy Chern class recovering the usual stringy Euler characteristic,
as explained above, and agreeing with the class introduced by deFernex
et al.

The $\omega$ flavor leads to a technical difficulty, which may make
the celestial integral (and hence stringy Chern classes) {\em
undefined\/} if the singularities of $X$ are not mild enough---the
technical condition is that they should be {\em log terminal.\/}
Whether stringy classes (or more generally celestial integrals) may be
defined for varieties with more general singularities is an open
question, see \S\ref{defin}.

One way out of this bind is to choose a different notion of relative
canonical divisor in the main set-up. For example, one can {\em
avoid\/} taking the double-dual, leading to the $\omega$ flavor as
mentioned above; this leads to a different notion (which I call the
$\Omega$ flavor of the integral), which is defined for arbitrarily
singular varieties. While this yields a stringy Chern class for
arbitrary varieties, the meaning of this class (for example vis-a-vis
the stringy Euler characteristic) has not been explored.


\section{Chern-Schwartz-MacPherson classes from celestial integrals}
\label{CSM}

\subsection{}
The stringy notion presented in \S\ref{stringy} amounts to taking the
identity manifestation of the integral of $0$ over the whole
modification system $\cC_X$ of the variety $X$. By the normalization
property ((\ref{norm}) in Theorem~\ref{props}), this yields the usual
Chern class of the tangent bundle of $X$ when $X$ is nonsingular.

There is a different natural way to use the same tool and define a
class generalizing $c(TX)\cap [X]$: embed $X$ into an ambient
nonsingular variety $M$, then compute the identity manifestation of
the celestial integral of $0$ {\em over the constructible subset\/}
$\cX$ determined by $X$:
$$\left(\integ[M]{0}{\cX}\right)_{\id}\quad.$$ 
With due care, this class can be defined in $A_*X$ (our definition of
the celestial integral would only place it in $A_*M$); remarkably, as
such it does {\em not\/} depend on the ambient variety~$M$. In fact:
$$\left(\integ[M]{0}{\cX}\right)_{\id}\overset!= \csm(X)\quad,$$
the {\em Chern-Schwartz-MacPherson\/} class of $X$.

This is a famous notion, going back to Marie-H\'el\`ene Schwartz
(\cite{MR35:3707}) and Robert MacPherson (\cite{MR50:13587}). In
MacPherson's construction (as recalled, for example, in
\cite{MR85k:14004}, \S19.1.7), one obtains in fact a natural
transformation $c_*$ from the functor of constructible functions (with
proper push-forward defined by Euler characteristic of fibers) to the
Chow group functor; applying $c_*$ to the constant function $\one_X$
defines the class $\csm(X)$.

\subsection{}
The connection between celestial integrals and Chern-Schwartz-MacPherson
classes mentioned above goes in fact much deeper. Given any divisor $\cD$
and any constructible subset $\cS$ of a modification system $\cC_X$, 
one may define a constructible function $I_X(\cD,\cS)$ by
$$I_X(\cD,\cS)(p):=\deg \left(\integ{\cD}{\cS\cap p}
\right)\quad;$$
here, $\cS\cap p$ is the constructible subset of $\cC_X$ obtained by
intersecting $\cS$ with inverse images of~$p$ through the system.

\begin{theorem}
$$\left(\integ{\cD}{\cS}\right)_{\id}=c_*(I_X(\cD,\cS))\quad.$$
\end{theorem}

Thus, celestial integrals and Chern-Schwartz-MacPherson classes are,
in a sense, equivalent information: each can be obtained from the
other.

Classes such as the stringy Chern class considered in \S\ref{stringy}
correspond, via MacPherson's natural transformation, to specific
constructible functions. These `stringy' characteristic functions
deserve much further study.

The apparatus of Chern-Schwartz-MacPherson classes is an important
ingredient in the construction in \cite{math.AG/0407314}.

\subsection{}
It should be noted that the definition of celestial integration (which
I will finally summarize in \S\ref{defin}) does {\em not\/} rely on
Chern-Schwartz-MacPherson classes; the latter are an honest subproduct
of the former. Thus, one could try to recover the main defining
features of Chern-Schwartz-MacPherson classes from celestial
properties.

For example, I would like to venture the guess that {\em the
covariance property of Chern-Schwartz-MacPherson classes is a facet of
the change-of-variable formula for celestial integrals.\/} This should
mean that the change-of-variable formula is a Riemann-Roch theorem in
disguise. As things stand now I don't even have a precise version of
this `guess' to offer, and I will have to leave it at the stage of
half-baked speculations.


\section{The definition}\label{defin}

A summary of celestial integration {\em without\/} a definition of 
this notion would be incomplete, even though I have tried to 
defend the idea that the definition itself is less important
than the fact alone that such a notion exists---in practice, the
normalization and change-of-variable properties suffice for
interesting applications and do not require the actual definition of
the integral to be appreciated.

In any case, here is the definition. Given a divisor $\cD$ and a
constructible (say closed and proper, for simplicity) subset $\cS$ of 
the modification system $\cC_X$, embedded resolution of singularities
ensures that there is an object $\pi: V_\pi \to X$ in $\cC_X$,
a normal crossing divisors $E$ with nonsingular components
$E_j$, $j\in J$, and divisors $D_\pi$, $S_\pi=\cup_{j\in J_S}E_j$
of $V_\pi$ such that:
\begin{itemize}
\item $\cD$ is represented by $(\pi, D_\pi)$;
\item $\cS$ is represented by $(\pi, S_\pi)$;
\item $D_\pi+K_\pi=\sum_{j\in J} m_j E_j$, with $m_j\in \Qbb$.
\end{itemize}
This set of data depends on the chosen notion of relative canonical 
divisor $K_\pi$. {\em Assume that all coefficients $m_j$ are $>-1$.}

\begin{defin}\label{expr}
$$\left(\integ{\cD}{\cS}\right)_\pi:=c(\Omega_{V_\pi}(\log E)^\vee)
\cap \sum_{I\subset J,\,I\cap J_S\ne \emptyset} 
\frac{[\cap_{i\in I} E_i]}{\prod_{i\in I}(1+m_i)}$$
\end{defin}

This expression defines the manifestation of the integral on all varieties
such as $V_\pi$, in which the data $\cD$, $\cS$ are `resolved' by a 
normal-crossing divisor. The manifestation on any other variety is
obtained by push-forward, compatibly with the requirement that the
celestial integral is an element of the inverse limit $A_*\cC_X$.

The obvious difficulty with this definition is that it is not at all
clear that it should not depend on the chosen $\pi$ used to resolve
the given data.  In motivic integration, similar expressions are
obtained {\em a posteriori,\/} and compute intrinsically defined
objects, hence it is clear that they do not depend on the choices. In
celestial integration I have to prove the necessary independence
explicitly, directly from Definition~\ref{expr}.

\begin{theorem}\label{limit}
If all $m_j$ are $>-1$, then the given expression does define an
element of the inverse limit $A_*\cC_X$.
\end{theorem}

This is proved by applying the factorization theorem of
\cite{MR2003c:14016}, which reduces this claim to a computation across
blow-ups along nonsingular centers. Manipulating the expressions is a
somewhat messy, but manageable, exercise in standard intersection
theory.

The independence requires that all $m_j>-1$ (even though the
expression in Definition~\ref{expr} makes sense as soon as no $m_j$ is
$=-1$); this is where the singularities of $X$ may play a r\^ole for
the particular case $\cD=0$, as I discussed in \S\ref{diffi}: the
restriction $m_j>-1$ in this case amounts to the requirement that $X$
be log~terminal.

The difficulty arising if some $m_j\le -1$ is that in a chain of
varieties connecting two varieties where the data is resolved, one may
appear for which the expression in Definition~\ref{expr} does not make
sense, for the mundane reason that one of the denominators in the
expression may vanish.

This problem arises in many different contexts, of which celestial
integration is but one instance (see for example
\cite{math.AG/0401374}, \S8, Question~I).  While there is a feeling
that the obstacle is technical rather than conceptual, it has opposed
stubborn resistance to the attempts made so far to overcome it, 
and examples such as the one presented in \S3.4 in \cite{Veys} 
suggest that the issue may be more fundamental than initially
expected.

The question of exactly which celestial integrals are well-defined
outside the range specified in Theorem~\ref{limit} is subtle and 
difficult. Answering this question is a worthwhile challenge: the 
present state of affairs limits the scope of the definition of certain 
key celestial integrals and hence, as pointed out in \S\ref{applic}
and~\S\ref{stringin}, of some potentially interesting applications.




\begin{thebibliography}{AKMW02}

\bibitem[AKMW02]{MR2003c:14016}
Dan Abramovich, Kalle Karu, Kenji Matsuki, and Jaros{\l}aw W{\l}odarczyk.
\newblock Torification and factorization of birational maps.
\newblock {\em J. Amer. Math. Soc.}, 15(3):531--572, 2002.

\bibitem[Alu]{math.AG/0407150}
Paolo Aluffi.
\newblock {Modification systems and integration in their Chow groups}.
\newblock FSU04-10, arXiv:math.AG/0407150.

\bibitem[Cra]{math.AG/9911179}
Alastair Craw.
\newblock {An introduction to motivic integration}, arXiv:math.AG/9911179.
\newblock {T}o appear in `{S}trings and {G}eometry'; M. Gross, J. Gauntlett and
  M. Douglas, eds.

\bibitem[dFLNU]{math.AG/0407314}
Tommaso de~Fernex, Ernesto Lupercio, Thomas Nevins, and Bernardo Uribe.
\newblock {Stringy Chern classes of singular varieties}, arXiv:math.AG/0407314.

\bibitem[DL01]{MR1905328}
Jan Denef and Fran{\c{c}}ois Loeser.
\newblock Geometry on arc spaces of algebraic varieties.
\newblock In {\em European Congress of Mathematics, Vol. I (Barcelona, 2000)},
  volume 201 of {\em Progr. Math.}, pages 327--348. Birkh\"auser, Basel, 2001.

\bibitem[Ful84]{MR85k:14004}
William Fulton.
\newblock {\em Intersection theory}.
\newblock Springer-Verlag, Berlin, 1984.

\bibitem[Loo02]{MR2003k:14010}
Eduard Looijenga.
\newblock Motivic measures.
\newblock {\em Ast\'erisque}, (276):267--297, 2002.
\newblock S\'eminaire Bourbaki, Vol.\ 1999/2000.

\bibitem[Mac74]{MR50:13587}
Robert~D. MacPherson.
\newblock Chern classes for singular algebraic varieties.
\newblock {\em Ann. of Math. (2)}, 100:423--432, 1974.

\bibitem[Sch65]{MR35:3707}
Marie-H{\'e}l{\`e}ne Schwartz.
\newblock Classes caract\'eristiques d\'efinies par une stratification d'une
  vari\'et\'e analytique complexe. {I}.
\newblock {\em C. R. Acad. Sci. Paris}, 260:3262--3264, 1965.

\bibitem[Tot00]{MR2001g:58037}
Burt Totaro.
\newblock Chern numbers for singular varieties and elliptic homology.
\newblock {\em Ann. of Math. (2)}, 151(2):757--791, 2000.

\bibitem[Veya]{math.AG/0401374}
Willem Veys.
\newblock {Arc spaces, motivic integration and stringy invariants},
  arXiv:math.AG/0401374.

\bibitem[Veyb]{Veys}
Willem Veys.
\newblock {On motivic principal value integrals}.
\newblock preprint.

\end{thebibliography}
\end{document}